\documentclass[final,11pt]{amsart}
\usepackage{ amsmath, amsthm, amsfonts, graphicx}
\usepackage[dvipsnames,usenames]{color}
\usepackage{verbatim}

\theoremstyle{definition}
\newtheorem{theorem}{Theorem}[section]
\newtheorem*{theorem*}{Theorem}
\newtheorem{pro}[theorem]{Proposition}
\newtheorem{Def}[theorem]{Definition}
\newtheorem{lem}[theorem]{Lemma}

\newtheorem{cor}[theorem]{Corollary}

\theoremstyle{definition}
\newtheorem*{Def*}{Definition}

\newtheorem{Rem}[theorem]{Remark}

\numberwithin{equation}{section}

\newcommand{\eus}{Euclidean space}

\newcommand{\gss}{graphical self-shrinker}

\newcommand{\mc}{mean curvature}
\newcommand{\mcf}{mean curvature flow}
\newcommand{\maxp}{maximum principle}
\newcommand{\ms}{minimal surface}
\newcommand{\nbd}{normal bundle}
\newcommand{\pf}{parallel form}
\newcommand{\pvgp}{polynomial volume growth property}

\newcommand{\ssk}{self-shrinker}

\newcommand{\sff}{second fundamental form}
\renewcommand{\rm}{\mathrm}

\newcommand{\tg}{totally geodesic}

\newcommand{\wrt}{with respect to}

\newcommand{\bpo}{\begin{pro}}
\newcommand{\epo}{\end{pro}}
\newcommand{\be}{\begin{equation}}
\newcommand{\ene}{\end{equation}}
\newcommand{\br}{\begin{Rem}}
\newcommand{\er}{\end{Rem}}
\newcommand{\bl}{\begin{lem}}
\newcommand{\el}{\end{lem}}
\newcommand{\bd}{\begin{Def}}
\newcommand{\ed}{\end{Def}}
\newcommand{\ben}{\begin{enumerate}}
\newcommand{\een}{\end{enumerate}}
\newcommand{\bp}{\begin{proof}}
\newcommand{\ep}{\end{proof}}
\newcommand{\beq}{\begin{equation*}}
\newcommand{\eeq}{\end{equation*}}
\newcommand{\bear}{\begin{eqnarray*}}
\newcommand{\eear}{\end{eqnarray*}}
\newcommand{\bt}{\begin{theorem}}
\newcommand{\et}{\end{theorem}}
\newcommand{\bst}{\begin{split}}
\newcommand{\est}{\end{split}}

\newcommand{\bal}{\begin{aligned}}
\newcommand{\eal}{\end{aligned}}

\renewcommand{\P}{\partial}
\newcommand{\F}[2]{\frac{#1}{#2}}
\newcommand{\la}{\langle}
\newcommand{\ra}{\rangle}

\newcommand{\R}{\mathbb{R}}

\newcommand{\lk}{\{}
\newcommand{\rk}{\}}

\newcommand{\bnb}{\bar{\nabla}}
\newcommand{\nb}{\nabla}

\newcommand{\Og}{\ast\Omega}

\newcommand{\Lp}{Laplacian}

\newcommand{\ex}{e^{-\F{|\vec{F}|^{2}}{4}}}
\newcommand{\cod}{co-dimension}
\newcommand{\Bt}{Bernstein}

\def\XXint#1#2#3{{\setbox0=\hbox{$#1{#2#3}{\int}$}
    \vcenter{\hbox{$#2#3$}}\kern-.5\wd0}}

\makeatletter
\def\@citestyle{\m@th\upshape\mdseries}
\def\citeform#1{{\bfseries#1}}
\def\@cite#1#2{{%
  \@citestyle[\citeform{#1}\if@tempswa, #2\fi]}}
\@ifundefined{cite }{%
  \expandafter\let\csname cite \endcsname\cite
  \edef\cite{\@nx\protect\@xp\@nx\csname cite \endcsname}%
}{}
\makeatother
\begin{document}
\title{A Bernstein type result for graphical self-shrinkers in $\R^{4}$}
\author{Hengyu Zhou}
\address[H. ~Z.]{Department of Mathematics, Sun Yat-sen University, 510275,
Guangzhou, P. R. China}
\email{hyuzhou84@yahoo.com}

\begin{abstract}
	%In this paper, we study nonparametric mean curvature type flows in $M\PLH\R$ which are represented as graphs $(x,u(x,t))$ over a domain in a Riemannian manifold $M$ with prescribed contact angle. The speed of $u$ is
	%the mean curvature speed minus a $\psi(x,u,Du)$ with an admissible condition. \\
	% \indent We mainly investigate two cases: $\psi(x,u, Du)\equiv 0$ with
	% vertical contact angle and $\psi(x,u,Du)=h(x,u)\sqrt{1+|Du|^2}$ with $h_u(x,u)\geq h_0>0$. The long time existence and uniformly convergence are established.
	%   In the first case the result is optimal since in appendix A we present a counterexample for the corresponding long time existence in a warped product manifold. In the second case the uniform convergence gives the solution to Capillary problems with positive gravity. We also obtain a gradient estimate of mean curvature type equations with Neumann boundary condition provided that $\psi(x,u,Du)$ is admissible. \\
	% \indent At last we treat a graph over a convex domain in a Riemannian surface $M^2$ with nonnegative Ricci curvature. With an assumption about the curvature of the boundary its nonparametric mean curvature flow in $M^2\PLH \R$ will approach uniformly to a translating solution to mean curvature flow as $t$ goes to $\infty$.
	Self-shrinkers are important geometric objects in the evolution of mean curvature flows, while the Bernstein Theorem is one of the most profound results in minimal surface theory. We prove a Bernstein type result for {\gss} surfaces with {\cod} two in $\R^4$. Namely under certain natural conditions on the Jacobian of any smooth map from $\R^2$ to $\R^2$ we show that the self-shrinker which is the graph of this map must be affine linear through 0. The proof relies on the derivation of structure equations of graphical {\ssk}s in terms of the {\pf} and the existence of some positive functions on
	self-shrinkers related to these Jacobian conditions.
\end{abstract}

\maketitle

\section{Introduction}
A smooth submanifold $\Sigma^n$ in $\R^{n+k}$ is a {\it {\ssk}} if the equation
\be\label{eq:stru}
\vec{H}+\F{1}{2}\vec{F}^{\bot}=0
\ene
holds for any point vector $\vec{F}$ on $\Sigma^n$. Here $\vec{H}$ is the {\mc} vector of $\Sigma^{n}$ and $\bot$ is the projection of
$\vec{F}$ into the {\nbd} of $\Sigma^{n}$.\\
\indent
Self-shrinkers are important in the study of {\mcf}s for at least two reasons. First if $\Sigma$ is a {\ssk}, it is easily checked that
\beq
\Sigma_{t}=\sqrt{-t}\Sigma,\quad  -\infty < t  <0
\eeq
is a solution to the {\mcf}. Hence {\ssk}s are self-similar solutions to
the {\mcf}. On the other hand by Huisken (\cite{Hui90}) the blow-ups around a type I singularity converge weakly to nontrivial {\ssk}s after rescaling and choosing subsequences. Because of the parabolic maximum principle the finite time singularity of {\mcf}s for initial compact hypersurfaces is unavoidable. Therefore it is desirable to classify {\ssk}s under various geometric conditions.
%%%%%%%%%%%%
\subsection{Motivation}
The rigidity of graphical minimal submanifolds in Euclidean spaces is summarized as the {\Bt} theorem. In this subsection we always assume that $f$ is a smooth map from $\R^{n}$ into $\R^{k}$ and $\Sigma =(x,f(x))$ is the graph of $f$. Let $Df$ denote the gradient of $f$. The {\Bt} theorem states that if $\Sigma$ is minimal, then $\Sigma$ is totally geodesic under the following conditions:
\begin{enumerate}
	\item  $n\leq 7$ and $k=1$ by \cite{SJ68};
	\item  any $n$ and $k=1$ with $|Df|=o(\sqrt{|x|^2+|f|^2})$ as $|x|\rightarrow \infty$ by \cite{EH90};
	\item  $n=2$ and any $k$ with $|Df|\leq C$  by \cite{CO67};
	\item $n=3$ and any $k$ with $|Df|\leq C$ by \cite{FD80};
	\item any $n\geq 2$ and $k\geq 2$ with more restrictive conditions (See Remark \ref{rm:minimal}) by \cite{WMT03}( also see \cite{JX99}, \cite{JXY13}).
\end{enumerate}
On the other hand a self-shrinker is minimal in weighted Euclidean space $(\R^{n+k}, e^{-\F{|x|^{2}}{2n}}dx^{2})$ where $dx^{2}$ is the standard Euclidean metric (\cite{CM12}). The interests in Bernstein type results for graphical self-shrinkers are revived due to the works of Ecker-Huisken (\cite{EH89}) and Wang (\cite{WL11}). They showed that a graphical {\ssk} $\Sigma$ with {\cod} one is a hyperplane through $0$ without any restriction on dimension $n$. Combining this with the historical results above, naturally we are interested in the rigidity of \emph{graphical {\ssk} surfaces} with {\cod} $k\geq 2$.\\
\indent There are two main difficulties to study graphical self-shrinker surfaces with higher codimension. First the techniques for minimal submanifolds in Euclidean space are generally not available. In the case of {\ssk} surfaces, there are no corresponding harmonic functions (\cite{CO67}) and monotonicity formula of the tangent cone at infinity in the minimal surface theory (\cite{SJ68}, \cite{FD80}, also \S 17 in \cite{SL83}). Second the contrast
between the hypersurface and higher co-dimensional submanifolds are another obstacle to study self-shrinkers with higher {\cod}. In the hypersurface case the normal bundle is trivial and the mean curvature is a scalar function. In the higher {\cod}
case the normal bundle can be highly non-trivial. In general the computations related to mean curvature and second fundamental form in this situation are very involved except few cases. \\
\indent However recent progresses on {\ssk}s and graphical {\mcf}s provide new tools to overcome these obstacles under some conditions. By \cite{CZ13} and \cite{DX13} there is a {\pvgp} for completely immersed, \emph{proper} {\ssk}s (Definition \ref{pvgp}). With this property, the integration technique gives good estimates if there are well-behaved structure equations satisfied by {\ssk}s (Lemma \ref{pro:est}). On the other hand graphical {\ssk}s have a lot of structure equations in terms of parallel form (Theorem \ref{thm:sk:eq} and \ref{pro:Laplacian}). This approach is inspired by the works of \cite{WMT01, WMT02} and \cite{TW04} to investigate graphical {\mcf}s with arbitrary codimension in product manifolds. They obtained evolution equations of the Hodge star of parallel forms along {\mcf}s (see Remark \ref{rm:star}). \\
\indent The contribution in this paper is to apply the parallel form's theory into the study of graphical self-shrinkers. In $\R^{4}$ both of codimension and dimension of a graphical self-shrinker surface are two. Then $\R^{4}$ provides four parallel 2-forms to reflect various properties of a graphical self-shrinker surface, which is explained in the next subsection.
%%%%%%%%%%%%
\subsection{Statement of the main result}
Suppose $f=(f_{1}(x_1,x_2), f_{2}(x_1,x_2))$ is a smooth map from $\R^{2}$ into $\R^{2}$. Then its Jacobian $J_f$ is given by
$$J_{f}=(\F{\P f_{1}}{\P x_{1}}\F{\P f_{2}}{\P x_{2}}-\F{\P f_{1}}{\P x_{2}}\F{\P f_{2}}{\P x_{1}})$$ The main result of this note is given as follows.
\bt\label{main}
Suppose $f:\R^2\rightarrow \R^2$ is a smooth map with its Jacobian $J_f$ satisfying (1): $J_f >-1$ for all $x$ or (2): $J_f<1$ for all $x$. If its graph is a self-shrinker in $\R^4$, then its graph is a two dimensional plane through 0.
\et
\br Notice that in our setting the {\cod} of $\Sigma$ is two.
\er
\br In \cite{DW10}, the authors studied the graphical self-shrinker of $f(x)$ in $\R^{n+m}$ where $f:\R^n\rightarrow \R^m$. While their approach is very promising in arbitrary codimension, it requires that the eigenvalues $\{\lambda_k\}$ of the maps $f$ satisfy $|\lambda_i\lambda_j|\leq 1$ for $i\neq j$. Notice that in Theorem \ref{main} we have $|J_f|=|\lambda_1\lambda_2|$ and $n=m=2$.\\
\indent Our paper relaxes the conditions on eigenvalues for graphical self-shrinkers, since we are able to use special geometry of 4-dimensions to prove our main results. Such geometry includes the existence of two parallel forms $dx_1\wedge dx_2$ and $dx_3\wedge dx_4$. See \S 3.1.\\
\indent It seems very difficult to generalize our approach to treat ever higher codimensional cases without further geometric inputs.
\er
A special case for Theorem \ref{main} is that if $f$ is a diffeomorphism on $\R^2$ the {\gss} of $f$ in $\R^{4}$ is totally geodesic. One can compare this theorem with the results (\cite{HSV09, HSV11}). Those authors obtained the rigidity of graphical minimal surfaces in $\R^{4}$ assuming bounded Jacobian.\\
\indent Let us explain conditions (1) and (2) in more details. Let $\Sigma$ be the {\gss} in Theorem ~\ref{main}. We
take $(x_{1},x_{2}, x_{3}, x_{4})$ as the coordinate of $\R^{4}=\R^{2}\times \R^{2}$. Let
$\eta_{1}=dx_{1}\wedge dx_{2}$, $\eta_{2}=dx_3\wedge dx_4$, $\eta' = \eta_1+\eta_2$ and
$\eta''=\eta_1 -\eta_2$. First we choose a proper orientation on $\Sigma$
such that $\ast\eta_{1}>0$ (Def. \ref{Def:star}). The direct computation shows that $\ast\eta'=(1+J_{f})\ast\eta_{1}$ and $\ast\eta''=(1-J_{f})\ast\eta_{1}$
(Lemma \ref{lm:hodge-star}). Then conditions (1) and (2) correspond to $\ast\eta'>0$ and $\ast\eta''>0$ respectively. When
{\it both} conditions (1) and (2) are satisfied, then we have $|J_f| < 1$ which means the map $f$ is area-decreasing. In \cite{WMT02, TW04}, assuming $f$ is area-decreasing, together with additional curvature conditions, they showed the {\mcf} of the graph of some smooth map stays
graphical and exists for all time.

Note that the usual {\maxp} does not apply to our non-compact submanifold. Our main technical tool to treat
this problem is Lemma ~\ref{pro:est} where we use a cutoff function and apply the Divergence Theorem:
a technique also used in \cite{WL11} for the case of hypersurface. The crucial condition is the {\pvgp} of {\gss}s. It can be of independent interest and 
we state it in a more general formulation (Lemma \ref{pro:est}).
%%%%%%%%%%%%
\subsection{Plan of the paper}
In \S2 we discuss the {\pf} and the geometry of {\gss}s. The structure equation of {\ssk}s in terms of
{\pf}s is summarized in Theorem \ref{thm:sk:eq}. In \S3 we apply Theorem \ref{thm:sk:eq} to the cases
of $*\eta'$ and $*\eta''$. For example in Theorem \ref{pro:Laplacian} we derive that
$*\eta' =\eta' (e_{1}, e_{2})$ satisfies the equation
\be\label{eq:key}
\Delta(*\eta') + *\eta'((h^{3}_{1k}- h^{4}_{2k})^{2}+(h^{4}_{1k}+h^{3}_{2k})^{2})-\F{1}{2}\la\vec{F},
\nb(*\eta') \ra =0,
\ene
where $h^{\alpha}_{ij}$ are the {\sff} and $\Delta$ ($\nb$) is the Laplacian (covariant derivative) of $\Sigma$.
With the {\pvgp} Lemma ~\ref{pro:est} implies that
\beq
((h^{3}_{1k}- h^{4}_{2k})^{2}+(h^{4}_{1k}+h^{3}_{2k})^{2})\equiv 0,
\eeq
if $*\eta'$ is a positive function. This implies that the {\gss} is minimal. Similar conclusion can be achieved for
$*\eta''$. We then show that it is actually {\tg}.
%%%%%%%%%%%%

%%%%%%%%%%%%%%%%%%%%%%%%%%%%
\section{Parallel Forms}\label{se:paral}
The {\pf}s in Euclidean space play the fundamental role in this paper. We will record many structure equations
of {\ssk}s in terms of {\pf}s. These equations can be quite general. We will present these results for submanifolds of arbitrary (co-)dimensions in general Riemannian manifolds in \S2.1 and then restrict to {\ssk}s in
Euclidean spaces in \S2.2.
%%%%%%%%%%%%%
\subsection{Parallel forms and their Hodge star}
We will adapt the notation in \cite{WMT08}, \cite{LL11}. Assume that $N^n$ is a smooth
$n$-dimensional submanifold in a Riemannian manifold $M^{n+k}$ of dimension $n+k$.  We denote an orthonormal basis of the tangent bundle of $N$ by $\lk e_{i}\rk_{i=1}^{n}$  and denote an orthonormal basis of the normal bundle of $N$ by $\lk e_{\alpha}\rk_{\alpha=n+1}^{n+k}$. The Riemann
curvature tensor of $M$ is defined by
\beq
R(X, Y)Z=-\bnb_{X}\bnb_{Y}Z+\bnb_{Y}\bnb_{X}Z+\bnb_{[X, Y]}Z,
\eeq
for smooth vector fields $X, Y$ and $Z$. The {\sff} $A$ and the {\mc} vector $\vec{H}$ are defined as
\begin{align}
A(e_{i}, e_{j})&=(\bnb_{e_{i}}e_{j})^{\bot}=h^{\alpha}_{ij}e_{\alpha}\\
\vec{H} &=(\bnb_{e_{i}}e_{i})^{\bot}=h^{\alpha}_{ii}e_{\alpha}=h^{\alpha}e_{\alpha}.
\end{align}
Here we used Einstein notation and $h^{\alpha}=h^{\alpha}_{ii}$.

Let $\nb$ be the covariant derivative of $\Sigma$ {\wrt} the induced metric. Then $\nb^{\bot} A$ can be written
as follows:
\be
\nb^{\bot}_{e_{k}}A(e_{i}, e_{j})=h^{\alpha}_{ij,k} e_{\alpha}.
\ene
Note that $h^{\alpha}_{ij,k}$ is not equal to $e_{k}(h^{\alpha}_{ij})$ unless $\Sigma$ is a hypersurface. In fact
we have
\bl  \label{lm:cd} $h^{\alpha}_{ij,k}$ takes the following form:
\be\label{eq:derf}
h^{\alpha}_{ij,k} = e_{k}(h^{\alpha}_{ij})+h^{\beta}_{ij}\la e_{\alpha}, \bnb_{e_{k}}e_{\beta}\ra
-C_{ki}^{l}h^{\alpha}_{lj}-C_{kj}^{l}h_{li}^{\alpha},
\ene
where $\nb_{e_{i}}e_{j}=C_{ij}^{k}e_{k}$.
\el
\bp
By its definition
\beq
h^{\alpha}_{ij,k}=\la \nb^{\bot}_{e_{k}}A(e_{i}, e_{j}), e_{\alpha}\ra.
\eeq
The conclusion follows from expanding $\nb^{\bot}_{e_{k}}A(e_{i},e_{j})$
\begin{align*}
h^{\alpha}_{ij,k}&=\la \bnb_{e_{k}}(A(e_{i},e_{j})), e_{\alpha}\ra
-\la A(\nb_{e_{k}}e_{i}, e_{j}), e_{\alpha}\ra -  \la A(e_{i},\nb_{e_{k}}e_{j}), e_{\alpha}\ra\\
&=\la \bnb_{e_{k}}(h^{\beta}_{ij}e_{\beta}), e_{\alpha}\ra
-C_{ki}^{l}h^{\alpha}_{lj}-  C_{kj}^{l}h_{li}^{\alpha}\\
&= e_{k}(h^{\alpha}_{ij})+h^{\beta}_{ij}\la e_{\alpha}, \bnb_{e_{k}}e_{\beta}\ra
-C_{ki}^{l}h^{\alpha}_{lj}-C_{kj}^{l}h_{li}^{\alpha}.
\end{align*}
\ep
For later calculation we recall that the
Codazzi equation is
\be\label{eq:cdz}
R_{\alpha ikj} = h^{\alpha}_{ij,k}-h^{\alpha}_{ik,j},
\ene
where $R_{\alpha ikj}=R(e_{\alpha}, e_{i}, e_{k}, e_{j})$.
\begin{Def}\label{Def:star} An
	$n$-form $\Omega$ is called {\it parallel} if $\bnb \Omega = 0$ where $\bnb $ is the covariant derivative
	of $M$. \\
	\indent The  Hodge star $\ast\Omega$ on $N$ is
	defined by
	\be\label{star}
	\ast\Omega=\F{\Omega(X_{1}, \cdots, X_{n})}{\sqrt{det(g_{ij})}}
	\ene
	where $\{X_1, \cdots, X_n\}$ is a local frame on $N$ and $g_{ij}=\la X_i, X_j\ra$.
\end{Def}
\br\label{rm:star}  We denote by $M$ the product manifold $N_1\times N_2$, $\Omega$ the volume form of $N_1$. Then $\Omega$ is a parallel form in $M$.
If $N$ is a graphical manifold over $N_1$, then $*\Omega >0$ on $N$ for an appropriate orientation. For example the graphical {\ssk} $\Sigma$ in \S 1.2 satisfies that $*\Omega >0$ on $\Sigma$ where $\Omega$ is $dx_1\wedge \cdots \wedge dx_n$.\\
\indent A crucial observation is that $\ast\Omega$ is independent of the frame
$\{X_{1},\cdots, X_{n}\}$ up to a fixed orientation. This fact greatly simplifies our calculation. When $\{X_1, \cdots, X_n\}$ is an orthonormal frame $\{e_{1},\cdots, e_{n}\}$, $\ast\Omega=\Omega(e_1, \cdots, e_n)$.\\
\indent  The evolution equation of $*\Omega$ along mean curvature flows is the key ingredient in (\cite{WMT02}).
\er
\br\label{rm:minimal} In \cite{WMT03} the author proved that suppose $\Sigma=(x,f(x))$ is minimal where $f:\R^{n}\rightarrow \R^{k}$ and there exists $0<\delta<1$ and $K>0$ such that $|\lambda_i\lambda_j|\leq 1-\delta$ and $*\Omega> K$, then $\Sigma$ is affine linear. \\
\indent Here $\{\lambda_i\}_{i=1}^{n}$ is the eigenvalue of $df$ and $\Omega$ is $dx_1\wedge\cdots\wedge dx_n$.
\er

The following equation \eqref{Lop} first appeared as equation (3.4) in \cite{WMT02} in the proof
of the evolution equation of $\ast\Omega$ along the {\mcf}. We provide a proof for the sake of
completeness.
\bpo \label{pro:1} Let $N^n$ be a smooth submanifold of $M^{n+k}$. Suppose $\Omega$ is a
parallel $n$-form and $R$ is the Riemann curvature tensor of $M$. Then
$\ast\Omega=\Omega(e_{1},\cdots, e_{n})$ satisfies the following equation:
\be\label{Lop}
\Delta(\ast\Omega)=-\sum_{i,k}(h^{\alpha}_{ik})^{2}\ast\Omega+\sum_{i}(h^{\alpha}_{,i}+
\sum_{k}R_{\alpha kik })\Omega_{i\alpha}
+2\sum_{i<j,k}h^{\alpha}_{ik}h^{\beta}_{jk}\Omega_{i\alpha,j\beta}.
\ene
Here $\Delta$ denotes the Laplacian on $N$ {\wrt} the induced metric, and
$h^{\alpha}_{,k}=h^{\alpha}_{ii,k}$. In the second group of terms, $\Omega_{i\alpha}=\Omega(\hat{e}_{1},\cdots, \hat{e}_{n})$ with $\hat{e}_{s}=e_{s}$ for $s\neq i$ and $\hat{e}_{s}=e_{\alpha}$ for $s=i$. In the last group of terms, $\Omega_{i\alpha, j\beta}=\Omega(\hat{e}_{1},\cdots, \hat{e}_{n})$ with $\hat{e}_{s}=e_{s}$ for $s\neq i,j$, $\hat{e}_{s}=e_{\alpha}$  for $s=i$ and $\hat{e}_{s}=e_{\beta}$ for $s=j$.
\epo
\bp
Recall that $\nb$ and $\bnb$ are the covariant derivatives of $N$ and $M$ respectively. Fix
a point $p$ on $\Sigma$ and assume that $\lk e_{1},\cdots, e_{n}\rk$ is normal at $p$ {\wrt} $\nb$.
Lemma ~\ref{lm:cd} implies that
\be\label{eq:drn}
\nb_{e_{i}}e_{j}(p)=0,\quad h^{\alpha}_{ij,k}(p)=e_{k}(h^{\alpha}_{ij})(p)
+h^{\beta}_{ij}\la e_{\alpha}, \bnb_{e_{k}}e_{\beta}\ra(p).
\ene
Since $\bnb\Omega = 0$, we have
\begin{align}
\nb_{e_{k}}(\ast\Omega)&=\Omega(\bnb_{e_{k}}e_{1}, \cdots, e_{n})+\cdots+\Omega(e_{1},\cdots,\bnb_{e_{k}}e_{n})\notag\\
&=\sum_{i}h^{\alpha}_{ik}\Omega_{i\alpha};\label{eq:grt}
\end{align}
For $\nb_{e_{k}}\nb_{e_{k}}(\ast\Omega)$ we get
\be\label{eq:sk2}
\nb_{e_{k}}\nb_{e_{k}}(\ast\Omega)=\sum_{i}e_{k}(h^{\alpha}_{ik})\Omega_{i\alpha}
+\sum_{i}h^{\alpha}_{ik}e_{k}(\Omega_{i\alpha}).
\ene
The second term in \eqref{eq:sk2} can be computed as
\begin{align*}
\sum_{i}h^{\alpha}_{ik}e_{k}(\Omega_{i\alpha})&=\sum_{i}  h^{\alpha}_{ik}\Omega(e_{1},\cdots, \bnb_{e_{k}}e_{\alpha},\cdots, e_{n})+2\sum_{i<j}h^{\alpha}_{ik}h^{\beta}_{jk}\Omega_{i\alpha, j\beta}\notag\\
&=\sum_{i,\alpha}-(h^{\alpha}_{ik})^{2}\ast\Omega+ h^{\beta}_{ ik}\la e_{\alpha}, \bnb_{e_{k}}e_{\beta}\ra \Omega_{i\alpha} +2\sum_{i<j}h^{\alpha}_{ik}h^{\beta}_{jk}\Omega_{i\alpha, j\beta}.
\end{align*}
Plugging this into \eqref{eq:sk2} yields that
\begin{align}
\nb_{e_{k}}\nb_{e_{k}}(\ast\Omega)
&=-\sum_{i,\alpha}(h^{\alpha}_{ik})^{2}\ast\Omega+2\sum_{i<j}h^{\alpha}_{ik}h^{\beta}_{jk}\Omega_{i\alpha,j\beta}+\sum_{i}h^{\alpha}_{ki,k}\Omega_{i\alpha}\notag\\
&=-\sum_{i,\alpha}(h^{\alpha}_{ik})^{2}\ast\Omega+2\sum_{i<j}h^{\alpha}_{ik}h^{\beta}_{jk}\Omega_{i\alpha,j\beta}+\sum_{i}(h^{\alpha}_{kk,i}+R_{\alpha kik})\Omega_{i\alpha}.\notag
\end{align}
In view of \eqref{eq:drn} and \eqref{eq:cdz} we can finally conclude that
\begin{align*}
\Delta(\ast\Omega(p))&=\nb_{e_{k}}\nb_{e_{k}}(\ast\Omega)(p)-\nb_{\nb_{e_{k}}e_{k}}(\ast\Omega)(p) \\
&=-\sum_{i,k,\alpha}(h^{\alpha}_{ik})^{2}\ast\Omega+2\sum_{i<j,k}h^{\alpha}_{ik}h^{\beta}_{jk}\Omega_{i\alpha, j\beta}+\sum_{i}(h^{\alpha}_{,i}+\sum_{k}R_{\alpha kik})\Omega_{i\alpha}.
\end{align*}
This is the conclusion.
\ep
%%%%%%%%%%%%%%%%%%%%%%%%%%
\subsection{Self-shrinkers in {\eus}}
In this subsection we only consider the case when $M^{n+k}$ is the Eucildean space and $N^n$ is a {\ssk}.
\bl \label{lm:slm} Let $\Omega$ be a parallel n-form in $\R^{n+k}$. Suppose $N^n$ is  an n-dimensional {\ssk}
in $\R^{n+k}$. Using the notation in Proposition \ref{pro:1} we have
\be
\sum_{i}  \Omega_{i\alpha} h^{\alpha}_{,i}=\F{1}{2}\la\vec{F}, \nb(\ast\Omega) \ra
\ene
where $\vec{F}$ is any point on $N^n$.
\el
\bp
As in \eqref{eq:drn} we assume that $\{e_{1},\cdots, e_{n}\}$ is normal at $p$. From \eqref{eq:grt}
we compute $\nb(\ast\Omega)$ as follows:
\beq
\nb (\ast\Omega) =\nb_{e_{k}}(\ast\Omega)e_{k}=(\sum_{i}h^{\alpha}_{ik}\Omega_{i\alpha})e_{k}.
\eeq
This leads to
\be \label{eq:sk3}
\F{1}{2}\la\vec{F},\nb(\ast\Omega) \ra=\F{1}{2}\la\vec{F},e_{k}\ra(\sum_{i}h^{\alpha}_{ki}\Omega_{i\alpha}).
\ene
Recall that $\vec{H} =h^{\alpha}e_{\alpha}$. Then $h^{\alpha}=-\F{1}{2}\la\vec{F}, e_{\alpha}\ra$ since
$\vec{H}+\F{1}{2}\vec{F}^{\bot}=0$.  Taking the derivative of $h^{\alpha}$ {\wrt} $e_{i}$ we get
\begin{align}
e_{i}(h^{\alpha})&=\F{1}{2}h^{\alpha}_{ik}\la\vec{F},e_{k}\ra-\F{1}{2}\la\vec{F},e_{\beta}\ra\la \bnb_{e_{i}}e_{\alpha}, e_{\beta}\ra\notag\\
&=\F{1}{2}h_{ik}^\alpha \la \vec{F}, e_k\ra+h^\beta \la \bnb_{e_i}e_\alpha, e_\beta\ra\notag\\
&=\F{1}{2}h^{\alpha}_{ik}\la\vec{F},e_{k}\ra-h^{\beta} \la \bnb_{e_{i}}e_{\beta}, e_{\alpha}\ra.\label{eq:cmp}
\end{align}
Here we applied $\la \bnb_{e_i}e_\alpha, e_\beta\ra=-\la \bnb_{e_i}e_{\beta},e_\alpha\ra$. 
Since we assume that $\nb_{e_{i}}e_{j}(p)=0$, \eqref{eq:derf} yields that $h^{\alpha}_{kk,i}(p)=e_{i}(h^{\alpha}_{kk})(p)+h^{\beta}_{kk} \la \bnb_{e_{i}}e_{\beta}, e_{\alpha}\ra(p)$. Then we conclude that
\beq
h^{\alpha}_{,i}(p)=e_{i}(h^{\alpha})(p)+h^{\beta} \la \bnb_{e_{i}}e_{\beta}, e_{\alpha}\ra(p).
\eeq
Comparing the above with \eqref{eq:cmp} we get $ h^{\alpha}_{,i}(p)=\F{1}{2}h^{\alpha}_{ik}\la\vec{F},e_{k}\ra(p)$.  The lemma follows from combining this with \eqref{eq:sk3}.
\ep
Using Proposition \ref{pro:1} and Lemma \ref{lm:slm} we obtain a series of structure equations of {\ssk}s
in terms of the {\pf}.
\bt\label{thm:sk:eq}(Structure Equation)
In $\R^{n+k}$ suppose $\Sigma$ is an $n$-dimensional {\ssk}. Let $\Omega$ be a parallel $n$-form, then
$\ast\Omega=\Omega(e_{1}, \cdots, e_{n})$ satisfies that
\begin{align}
\Delta (\Og) +(h^{\alpha}_{ik})^{2}\Og
-2\sum_{i< j} \Omega_{i\alpha,j\beta}h^{\alpha}_{ik}h^{\beta}_{jk}-
\F{1}{2}\la\vec{F},\nb(\ast\Omega\ra)=0,\label{eq:sk}
\end{align}
where $\vec{F}$ is the coordinate of the point on $\Sigma$ and  $\Omega_{i\alpha, j\beta}=\Omega(\hat{e}_{1},\cdots, \hat{e}_{n})$ with $\hat{e}_{s}=e_{s}$ for $s\neq i,j$, $\hat{e}_{s}=e_{\alpha}$  for $s=i$ and $\hat{e}_{s}=e_{\beta}$ for $s=j$.
\et
This theorem enables us to obtain various information of {\ssk}s for different {\pf}s. We will apply this idea
to our particular situation in the next section.
%%%%%%%%%%%%%%%%%%%%%%%%%%%
\section{Graphical {\ssk}s in ${\R^4}$}
From this section on we will focus on the {\gss}s in {\eus}. The structure equations of
graphical {\ssk}s will be derived in Theorem  \ref{pro:Laplacian}. The {\pvgp} plays an essential role in Lemma
~\ref{pro:est}, which is our main technical tool.
%%%%%%%%%%%%%%%%%%%%%%
\subsection{Structure equations for {\gss}s in $\R^{4}$} We consider the following four different
parallel 2-forms in $\R^4$:
\begin{align}
&\eta_{1}=dx_{1}\wedge dx_{2},  &\eta'=dx_{1}\wedge dx_{2}+dx_{3}\wedge dx_{4}\notag\\
&\eta_{2}=dx_{3}\wedge dx_{4} &\eta''=dx_{1}\wedge dx_{2}-dx_{3}\wedge dx_{4}\label{def:coor}
\end{align}
Recall that for a smooth map $f=(f_{1}(x_{1}, x_{2}),f_{2}(x_{1}, x_{2}))$ its Jacobian $J_{f}$ is
\be
J_{f}=(\F{\P f_{1}}{\P x_{1}}\F{\P f_{2}}{\P x_{2}}-\F{\P f_{1}}{\P x_{2}}\F{\P f_{2}}{\P x_{1}}) ;
\ene
\bl \label{lm:hodge-star} Suppose $\Sigma =(x, f(x))$ where $f:\R^{2}\rightarrow \R^{2}$ is a smooth map. Then on $\Sigma$ it holds that 
$$
*\eta_{2}= J_{f}*\eta_{1};
$$
\el
\bp Notice that $*\eta_1$ and $*\eta_2$ are independent of the choice of the local frame.  Denote by $e_1= \F{\P }{\P x_1}+\F{\P f_1}{\P x_1}\F{\P}{\P x_3}+\F{\P f_2}{\P x_1}\F{\P}{\P x_4}$,  $e_2= \F{\P }{\P x_2}+\F{\P f_1}{\P x_2}\F{\P}{\P x_3}+\F{\P f_2}{\P x_2}\F{\P}{\P x_4}$ and $g_{ij}=\la e_i, e_j \ra$. Then
\begin{align}
*\eta_2&=\F{dx_3\wedge dx_4(e_1, e_2)}{\sqrt{det(g_{ij})}}\notag\\
&=\F{J_f}{\sqrt{det(g_{ij})}}\notag\\
&=J_f*\eta_1; \notag
\end{align}
\ep
The above lemma is not enough to explore structure equations in Theorem \ref{thm:sk:eq}.
We need further information about the microstructure of a point on $\Sigma$.
\bl  Assume $f:\R^{2}\rightarrow \R^{2}$ is a smooth map. Denote by $df$ the differential of $f$. Then for any point $x$
\begin{enumerate}
	\item There exist oriented orthonormal bases $\{a_{1}, a_{2}\}$ and $\{a_{3}, a_{4}\}$ in $T_{x}\R^{2}$ and $T_{f(x)}\R^{2}$ respectively such that
	\be\label{eq:lin_map}
	df(a_{1})=\lambda_{1}a_{3},\qquad df(a_{2})=\lambda_{2}a_{4};
	\ene
	Here `oriented' means $dx_{i}\wedge dx_{i+1}(a_{i}, a_{i+1})=1$ for $i=1,3$ .
	\item Moreover we have $\lambda_{1}\lambda_{2}=J_{f}$.
\end{enumerate}
\el
\bp Fix a point $x$.  First we prove the existence of (1). By the Singular Value Decomposition Theorem (p.291 in \cite{ST07}) there exist two $2\times 2$ orthogonal matrices $Q_1, Q_2$ such that
$$
\begin{pmatrix}
\F{\P f_{1}}{\P x_{1}} & \F{\P f_{2}}{\P x_{1}}\\
\F{\P f_{1}}{\P x_{2}} & \F{\P f_{2}}{\P x_{2}}
\end{pmatrix}
= Q_1\begin{pmatrix}
\lambda'_{1} & 0 \\
0& \lambda'_{2}
\end{pmatrix} Q_2
$$
with $\lambda'_{1},\lambda'_2\geq 0$. \\
\indent Let $\lambda_{1}=\rm{det}(Q_1)\lambda'_{1}\rm{det}(Q_2)$, $\lambda_{2}=\lambda'_{2}$, $A =\rm{det}(Q_1)Q_1$ and $B=\rm{det}(Q_2)Q_2$. Thus $\rm{det}(A)=\rm{det}(B)=1$ and 
\be\label{eq:cord}
\begin{pmatrix}
	\F{\P f_{1}}{\P x_{1}} & \F{\P f_{2}}{\P x_{1}}\\
	\F{\P f_{1}}{\P x_{2}} & \F{\P f_{2}}{\P x_{2}}
\end{pmatrix}
= A \begin{pmatrix}
	\lambda_{1} & 0 \\
	0& \lambda_{2}
\end{pmatrix} B
\ene
We consider the new basis $(a_{1},a_{2})^{T} =A^T(\F{\P}{\P x_{1}}, \F{\P }{\P x_{2}})^{T}$, $(a_{3}, a_{4})^{T}=B(\F{\P}{\P x_{3}}, \F{\P }{\P x_{4}})^{T}$, then $dx_{1}\wedge dx_{2}(a_1, a_2)=1$ and $dx_{3}\wedge dx_{4}(a_3, a_4)=1$ ($A^T$ is the transpose of $A$). Moreover \eqref{eq:cord} implies that
$$
df(a_1, a_2)^{T}=\begin{pmatrix}
\lambda_{1} &0 \\
0 &\lambda_{2}
\end{pmatrix} (a_{3}, a_{4})^T
$$
Now we obtain (1).\\
\indent According to \eqref{eq:cord} we have $J_f=det(A)\lambda_1\lambda_2 det(B)=\lambda_1\lambda_2$. We arrive at (2). The proof is complete. 
\ep
\br The conclusion (2) does not depend on the special choice of $\{a_i\}_{i=1}^4$ which satisfies \eqref{eq:lin_map}. 
\er 
With these bases we construct the following local frame for later use.
\begin{Def}\label{def:cor}
	Fix a point $p=(x, f(x))$ on $\Sigma$. We construct a special orthonormal basis $\{e_{1},e_{2}\}$
	of the tangent bundle $T\Sigma$ and $\{e_{3}, e_{4}\}$ of the normal bundle $N\Sigma$ at follows. At the point $p$ we have for $i=1,2$:
	\begin{align}
	e_{i}=\F{1}{\sqrt{1+\lambda^{2}_{i}}}(a_{i}+\lambda_{i}a_{2+i});\quad
	e_{2+i}=\F{1}{\sqrt{1+\lambda^{2}_{i}}}(a_{2+i}-\lambda_{i}a_{i})\label{eq:cd};
	\end{align}
	where $\{a_{1}, a_2, a_3, a_4\}$ are from \eqref{eq:lin_map}.
\end{Def}
For a parallel 2-form $\Omega$ we have $\ast\Omega=\Omega(e_{1},e_{2})$. Applying \eqref{eq:cd}
and $\lambda_{1}\lambda_{2}=J_{f}$ direct computations show that $*\eta_{1},*\eta_{2}, *\eta'$ and
$*\eta''$ take the following form:
\begin{align}
*\eta_{1}&=\F{1}{\sqrt{(1+\lambda_{1}^{2})(1+\lambda_{2}^{2})}} > 0, \label{eta1}\\
*\eta_{2}&=\F{\lambda_{1}\lambda_{2}}{\sqrt{(1+\lambda_{1}^{2})(1+\lambda_{2}^{2})}},\\
*\eta'&=(1+J_{f})(*\eta_{1}),\label{eta-prime}\\
*\eta''&=(1-J_{f})(*\eta_{1}).\label{eq:ineta}
\end{align}
Now we have the
structure equations for {\gss}s in $\R^4$ as follows:
\bt\label{pro:Laplacian}
Suppose $f:\R^{2}\rightarrow \R^{2}$ is a smooth map and $\Sigma =(x, f(x))$ is a {\gss} in $\R^{4}$.
Using the notation in Definition ~\ref{def:cor} we have
\begin{align}
\Delta (*\eta_1)&+*\eta_{1}(h^{\alpha}_{ik})^{2}
-2*\eta_{2}(h^{3}_{1k}h^{4}_{2k}-h^{4}_{1k}h^{3}_{2k})-\F{1}{2}\la\vec{F},\nabla(*\eta_1)\ra=0;\label{eq:lp:1}\\
\Delta(*\eta_{2})&+*\eta_{2}(h^{\alpha}_{ik})^{2}-2*\eta_{1}(h^{3}_{1k}h^{4}_{2k}
-h^{4}_{1k}h^{3}_{2k})-\F{1}{2}\la\vec{F},\nabla(*\eta_{2})\ra=0;\label{eq:lp:2}\\
\Delta(*\eta')&+*\eta'((h^{3}_{1k}-h^{4}_{2k})^{2}+(h^{4}_{1k}
+h^{3}_{2k})^{2})-\F{1}{2}\la\vec{F},\nabla(*\eta')\ra=0;\label{eq:lp:3}\\
\Delta(*\eta'')&+*\eta''((h^{3}_{1k}+h^{4}_{2k})^{2}+(h^{4}_{1k}-h^{3}_{2k})^{2})
-\F{1}{2}\la\vec{F},\nabla(*\eta'')\ra=0,\label{eq:lp:4}
\end{align}
where $h^{\alpha}_{ij}=\la \bnb_{e_{i}}e_{j}, e_{\alpha}\ra$ are the second fundamental form of $\Sigma$, $\Delta$ and $\nabla$ are the {\Lp} and the {covariant derivative} of $\Sigma$ respectively.
\et
\bp First we consider the equation \eqref{eq:lp:1}. Applying the frame in \eqref{eq:cd}, the third term in
Theorem ~\ref{thm:sk:eq} becomes:
\begin{align*}
2(\eta_{1})_{i\alpha, j\beta}h^{\alpha}_{ik}h^{\beta}_{jk}&=2 dx_{1}\wedge
dx_{2}(e_{3}, e_{4})(h^{3}_{1k}h^{4}_{2k}-h^{3}_{2k}h^{4}_{1k})\\
&=2\F{\lambda_{1}\lambda_{2}}{\sqrt{(1+\lambda_{1}^{2})
		(1+\lambda_{2}^{2})}}(h^{3}_{1k}h^{4}_{2k}-h^{4}_{1k}h^{3}_{2k})\\
&=2*\eta_{2}(h^{3}_{1k}h^{4}_{2k}-h^{4}_{1k}h^{3}_{2k}).
\end{align*}
Here in the second line we used the fact that $dx_{1}\wedge dx_{2}(a_{1}, a_{2})=1$.
Plugging this into \eqref{eq:sk}, we obtain \eqref{eq:lp:1}. \\
\indent Similarly we obtain that
\begin{align}
2(\eta_{2})_{i\alpha, j\beta}h^{\alpha}_{ik}h^{\beta}_{jk}&=2 dx_{3}\wedge dx_{4}(e_{3}, e_{4})(h^{3}_{1k}h^{4}_{2k}-h^{3}_{2k}h^{4}_{1k})\notag\\
&=2\F{1}{\sqrt{(1+\lambda_{1}^{2})(1+\lambda_{2}^{2})}}(h^{3}_{1k}h^{4}_{2k}-h^{4}_{1k}h^{3}_{2k})\notag\\
&=2*\eta_{1}(h^{3}_{1k}h^{4}_{2k}-h^{4}_{1k}h^{3}_{2k}).\label{eq:if}
\end{align}
Here in the second line we used the fact that $dx_{3}\wedge dx_{4}(a_{3}, a_{4})=1$. Then
\eqref{eq:lp:2} follows from plugging \eqref{eq:if} into \eqref{eq:sk}.\\
\indent To show \eqref{eq:lp:3} we observe that $*\eta'=*\eta_1+*\eta_2$. Then plugging \eqref{eq:lp:1} into \eqref{eq:lp:2} we obtain that 
$$
\Delta (*\eta')+*\eta'\{\sum_{\alpha=3}^4\sum_{i,k=1}^2(h^\alpha_{ik})^2-2\sum_{k=1}^2(h^3_{1k}h^4_{2k}-h^4_{1k}h^3_{2k})\}-\F{1}{2}\la \vec{F},\nb(*\eta')\ra=0
$$ 
Thus \eqref{eq:lp:3} follows from the identity 
$$
\sum_{\alpha=3}^4\sum_{i,k=1}^2(h^\alpha_{ik})^2-2\sum_{k=1}^2(h^3_{1k}h^4_{2k}-h^4_{1k}h^3_{2k})=\sum_{k=1}^2(h^3_{1k}-h^4_{2k})^2+(h^4_{1k}+h^3_{2k})^2
$$
With a similar derivation we can show \eqref{eq:lp:4}. The proof is complete. 
\ep
%%%%%%%%%%%%%%%%%%%%%%%%%
\subsection{Volume growth for {\ssk}s}
We will state our main analytic tool in a more general setting since it may be of independent interest. In
this subsection we will consider {\gss}s of $n$-dimensional in $\R^{n+k}$.
\begin{Def}\label{pvgp}
	Let $N^n$ be a complete, immersed $n$-dimensional submanifold in $\R^{n+k}$. We say $N$ has
	the {\pvgp} if for any $r\geq 1$
	\beq
	\int_{N\cap B_{r}(0)} d vol \leq C r^n,
	\eeq
	where $B_{r}(0)$ is the ball in $\R^{n+k}$ centered at $0$ with radius $r$.
\end{Def}
Recently \cite{CZ13} and \cite{DX13}
showed the {\pvgp} is automatic under the following condition, but without the restriction of
dimension and codimension.
\bt (\cite{CZ13, DX13}) \label{thm:vgrowth}
If $N^n$ is a $n$-dimensional complete, \textit{properly} immersed {\ssk} in $\R^{n+k}$, then it
satisfies the {\pvgp}.
\et
\br \label{rm:vgth}
The properness assumption can not be removed. See Remark 4.1 in \cite{CZ13}.
\er
Notice that any {\gss} in {\eus} is embedded, complete and proper. Thus we have the following conclusion.
\begin{cor}
	Let $\Sigma =(x, f(x))$ be a smooth {\gss} in $\R^{4}$ where
	$f:\R^{2}\rightarrow \R^{2}$ is a smooth map. Then $\Sigma$ has the {\pvgp}.
\end{cor}
The following lemma is crucial for our argument:
\bl \label{pro:est}
Let $N^n \subset \R^{n+k}$ be a complete, immersed smooth $n$-dimensional submanifold
with at most polynomial volume growth. Suppose $g$ is a positive function and $K$ is a nonnegative function satisfying 
\be \label{eq:beq}
0\geq \Delta g -\F{1}{2}\la\vec{F}, \nb g\ra + Kg,
\ene
where $\Delta (\nb)$ is the Laplacian (covariant derivative) of
$N^n$ and $\vec{F}$ is the position vector of $N^n$. Then $g$ is a positive constant and $K\equiv 0$.
\el
\bp Fix $r\geq 1$. We denote by $\phi$ a compactly supported smooth function in $\R^{n+k}$
such that $\phi\equiv 1$ on $B_{r}(0)$ and $\phi\equiv 0$ outside of $B_{r+1}(0)$ with
$|\nb \phi|\leq |D\phi|\leq 2$. Here $D\phi$ and $\nb \phi$ are the gradient of $\phi$ in
$\R^{n+k}$ and $N^n$ respectively.

Since $g$ is positive, let $u=\log g$. Then the inequality \eqref{eq:beq} becomes
\beq
0\geq \Delta u -\F{1}{2} \la \vec{F}, \nb u\ra +(K+|\nb u|^{2}).
\eeq
Multiplying the righthand side of the above equation by $\phi\ex$ and integrating it on $N^n$ we get 
\begin{align}
0&\geq \int_{N}\phi^{2}div_{N}(\ex \nb u) + \int_{N}\phi^{2}\ex(K+|\nb u|^{2})\notag\\
&= -\int_{N} 2\phi\la\nb \phi, \nb u\ra \ex + \int_{N}\phi^{2}\ex(K+|\nb u|^{2})\notag\\
&\geq -\int_{N}2|\nb\phi|^{2}\ex +\int_{N}\phi^{2}\ex(K+\F{|\nb u|^{2}}{2})\label{eq:kt}
\end{align}
In \eqref{eq:kt} we used the inequality
\beq
|2\phi\la \nb \phi, \nb u\ra |\leq \F{\phi^{2}|\nb u|^{2}}{2}+2|\nb\phi|^{2}.
\eeq
Now we estimate, using the condition that $|\nb \phi|\leq |D\phi|\leq 2$,
\begin{align*}
\int_{N\cap B_{r}(0)}\ex(K+\F{|\nb u|^{2}}{2})&\leq  \int_{N}\phi^{2}\ex(K+\F{|\nb u|^{2}}{2})\\
&\leq \int_{N}2|\nb\phi|^{2}\ex \quad\text{ by \eqref{eq:kt}}\\
&\leq 8\int_{N\cap (B_{r+1}(0)\backslash B_{r}(0))}\ex\\
&\leq 8C (r+1)^{n}e^{-\F{r^{2}}{4}}.
\end{align*}
In the last line we use the fact that the submanifold $N^n$ has the {\pvgp}.

Letting $r$ go to infinity we obtain that
\beq
\int_{N}\ex(K +\frac{|\nb u|^2}{2})\leq 0.
\eeq
Since $K$ is nonnegative, we have $K\equiv \nb u \equiv 0$. Therefore $g$ is a positive constant.
\ep
%%%%%%%%%%%%%%%%%%%%%%%%%%%%
\subsection{The proof of Theorem \ref{main}}
Adapting to our case of {\gss} surfaces in $\R^4$ we are ready to prove Theorem \ref{main}.
\bp (of Theorem ~\ref{main})
We claim that $\Sigma$ is minimal under the assumptions. We prove this case by case.
\subsubsection*{Assuming Condition (1):} the equations \eqref{eta1} and \eqref{eta-prime} imply
that the {\pf} $\ast\eta'$ has the same sign as $1+ J_{f}$. Hence $\ast\eta'$ is a positive function. Moreover
from \eqref{eq:lp:3} $\ast\eta'$ satisfies
\beq
\Delta(\ast\eta')+(\ast\eta')((h^{3}_{1k}-h^{4}_{2k})^{2}+(h^{4}_{1k}+h^{3}_{2k})^{2})
-\F{1}{2}\la\vec{F},\nabla(\ast\eta')\ra=0.
\eeq
Here $\vec{F}=(x,f(x))$. Since $\Sigma$ has the {\pvgp}, using Lemma \ref{pro:est} we conclude that
\beq
(h^{3}_{1k}-h^{4}_{2k})^{2}+(h^{4}_{1k}+h^{3}_{2k})^{2}\equiv 0.
\eeq
We then obtain
\begin{align}
h^{3}_{11}&=h^{4}_{21}, \quad h^{3}_{22}=-h^{4}_{12}\label{eq:m1},\\
h^{4}_{11}&=-h^{3}_{21},\quad h^{4}_{22}=h^{3}_{12}\label{eq:m2}.
\end{align}
Then $\vec{H}=( h^{3}_{11}+h^{3}_{22})e_{3}+(h^{4}_{11}+h^{4}_{22})e_{4}\equiv 0$. So $\Sigma$ is a
{\ms}.
\subsubsection*{Assuming Condition (2):} this is similar to the above case. \eqref{eta1} and
\eqref{eq:ineta} imply that $\ast\eta''$ has the same sign as $1-J_{f}$. Thus $\ast\eta''$ is positive.
From \eqref{eq:lp:4} it also satisfies
\beq
\Delta(\ast\eta'')+(\ast\eta'')((h^{3}_{1k}+h^{4}_{2k})^{2}+(h^{4}_{1k}-h^{3}_{2k})^{2})
-\F{1}{2}\la\vec{F},\nabla(\ast\eta'')\ra=0
\eeq
Again we apply Lemma \ref{pro:est} to find that
\beq
(h^{3}_{1k}+h^{4}_{2k})^{2}+(h^{4}_{1k}-h^{3}_{2k})^{2}\equiv 0.
\eeq
Then we have
\begin{align}
h^{3}_{11}&=-h^{4}_{21}, \quad h^{3}_{22}=h^{4}_{12},\label{eq:m3}\\
h^{4}_{11}&=h^{3}_{21},\quad h^{4}_{22}=-h^{3}_{12}.\label{eq:m4}
\end{align}
Therefore we arrive at:
\beq
\vec{H}=( h^{3}_{11}+h^{3}_{22})e_{3}+(h^{4}_{11}+h^{4}_{22})e_{4}\equiv 0,
\eeq
which also means $\Sigma$ is minimal.\\
\indent Now $\Sigma$ is a {\gss} and minimal. From \eqref{eq:stru} we have $\vec{F}^{\bot}\equiv 0$
for any point $\vec{F}$ on $\Sigma$. For any normal unit vector $e_{\alpha}$ in the normal bundle of
$\Sigma$, we have
\be\label{F}
\la\vec{F}, e_{\alpha}\ra\equiv 0.
\ene
Take derivative {\wrt} $e_{i}$ for $i=1,2$ from \eqref{F} and we obtain
\begin{align*}
&\la\vec{F}, e_{1}\ra h^{\alpha}_{11}+\la\vec{F}, e_{2}\ra h^{\alpha}_{12}=0,\\
&\la\vec{F}, e_{1}\ra h^{\alpha}_{21}+\la\vec{F}, e_{2}\ra h^{\alpha}_{22}=0,
\end{align*}

Now assume $\vec{F}\neq 0$. Since $\vec{F}^{\bot}=0$, $(\la \vec{F}, e_{1}\ra, \la\vec{F}, e_{2}\ra)\neq (0, 0)$. According to the basic linear algebra we conclude that
\be\label{det}
h^{\alpha}_{11}h^{\alpha}_{22}-(h^{\alpha}_{12})^{2}=0
\ene
The minimality implies $h^{\alpha}_{11}=-h^{\alpha}_{22}$. Hence \eqref{det} becomes
$-(h^{\alpha}_{11})^{2}=(h^{\alpha}_{12})^{2}$. We find that $h^{\alpha}_{ij}=0$ for $i, j = 1,2$. Therefore $\Sigma$ is {\tg} except $\vec{F}= 0$.  Since $\Sigma$ is a graph, there is at most one point on $\Sigma$ such that $\vec{F} =0$. By the continuity of the {\sff} $\Sigma$ is {\tg} everywhere.\\
\indent Now $\Sigma$ is a plane. Provided $0$ is not on the plane, then we can find a point $\vec{F_0}$ in this plane which is nearest to $0$. It is easy to see that $\vec{F}_0=\vec{F}^{\bot}_0\neq 0$. This gives a contradiction because $\vec{F}^{\bot}=-\F{\vec{H}}{2}\equiv 0$. We complete the proof.
\ep
\section*{Acknowledgement}
This work is supported by the National Natural Science Foundation of China, No. 11271378 and No. 11521101. The author wishes to thank his thesis advisors, Prof. Zheng Huang and Prof. Yunping Jiang, for constant encouragements and valuable suggestions. He also wishes to thank Prof. Mao-Pei Tsui for insightful discussions. The author is very grateful to careful readings and helpful comments of the referees. 
\bibliographystyle{abbrv}	
\bibliography{SS_Ref}
 %\nocite{*}
 %\bibliography{ref_ms}

\end{document}